\documentclass[11pt]{article}
\title{
Sur la perfection de l'ensemble de Julia des
fonctions enti\`eres transcendantes
}

\font\sdopp=msbm9 scaled \magstep1
\def\CI {\sdopp {\hbox{C}}}
\author{Claudi Meneghin
}

\font\cir=wncyb10

\def\pe{{\cir\hbox{P}}}

\begin{document}
\maketitle
\bibliographystyle{plain} 
\parindent=8pt
\font\cir=wncyb10
\def\Iu{\cir\hbox{YU}}
\def\Ze{\cir\hbox{Z}}
\def\pe{\cir\hbox{P}}
\def\Ef{\cir\hbox{F}}
\def\CIRC{\mathop{\tt o}\limits}
\def\quan{\vrule height6pt width6pt depth0pt}
\def\QUAN{\ \quan}
\def\BETA{\mathop{\beta}\limits}
\def\GAMMA{\mathop{\gamma}\limits}
\def\VI{\mathop{v}\limits}
\def\UI{\mathop{u}\limits}
\def\VII{\mathop{V}\limits}
\def\WI{\mathop{w}\limits}
\def\ZETA{\mathop{Z}\limits}

\newtheorem{definition}{Definition}
\newtheorem{defi}[definition]{D\'efinition}
\newtheorem{lemma}[definition]{Lemma}
\newtheorem{lemme}[definition]{Lemme}
\newtheorem{proposition}[definition]{Proposition}
\newtheorem{theorem}[definition]{Theorem}        
\newtheorem{theoreme}[definition]{Th\'eor\`eme}        
\newtheorem{corollary}[definition]{Corollary}  
\newtheorem{corollaire}[definition]{Corollaire}  
\newtheorem{remark}[definition]{Remark}  
\newtheorem{remarque}[definition]{Remarque}
  
\font\sdop=msbm8
\def\ER {\sdopp {\hbox{R}}}
\def\ERp {\sdop {\hbox{R}}}
\def\QU {\sdopp {\hbox{Q}}}
\def\CI {\sdopp {\hbox{C}}}
\def\DI {\sdopp {\hbox{D}}}
\def\TI {\sdopp {\hbox{T}}}
\def\TU {\sdop {\hbox{T}}}
\def\EN{\sdopp {\hbox{N}}}
\def\ZETA{\sdopp {\hbox{Z}}}
\def\ES{\sdopp {\hbox{S}}}
\def\PI {\sdopp {\hbox{P}}}
\def\M{\hbox{\tt\large M}}
\def\N{\hbox{\boldmath{}$N$\unboldmath}} 
\def\P{\hbox{\boldmath{}$P$\unboldmath}} 
\def\tr{\hbox{\boldmath{}$tr$\unboldmath}} 
\def\f{\hbox{\large\tt f}} 
\def\F{\hbox{\boldmath{}$F$\unboldmath}} 
\def\G{\hbox{\boldmath{}$G$\unboldmath}} 
\def\L{\hbox{\boldmath{}$L$\unboldmath}} 
\def\h{\hbox{\boldmath{}$h$\unboldmath}} 
\def\e{\hbox{\boldmath{}$e$\unboldmath}} 
\def\g{\hbox{\boldmath{}$g$\unboldmath}} 
\def\u{\hbox{\boldmath{}$u$\unboldmath}} 
\def\v{\hbox{\boldmath{}$v$\unboldmath}} 
\def\U{\hbox{\boldmath{}$U$\unboldmath}} 
\def\V{\hbox{\boldmath{}$V$\unboldmath}} 
\def\W{\hbox{\boldmath{}$W$\unboldmath}} 
\def\id{\hbox{\boldmath{}$id$\unboldmath}} 
\def\alph{\hbox{\boldmath{}$\alpha$\unboldmath}} 
\def\bet{\hbox{\boldmath{}$\beta$\unboldmath}} 
\def\gam{\hbox{\boldmath{}$\gamma$\unboldmath}} 
\def\pphi{\hbox{\boldmath{}$\varphi$\unboldmath}} 
\def\ppsi{\hbox{\boldmath{}$\psi$\unboldmath}} 
\def\Ppsi{\hbox{\boldmath{}$\Psi$\unboldmath}} 
\def\brevve{}
\def\labelle #1{\label{#1}}
\def\quadras #1{ \hbox{\bf[}{#1}
\hbox{\bf]}}


Le but de cette note-ci est de prouver de fa\c con \'el\'ementaire
que l'ensemble de Julia 
d'une fonction enti\`ere transcendente \`a une variable 
complexe n'est pas vide et est parfait.
Ces r\'esultats sont bien s\^ur d\'eja connus: la
premi\`ere preuve remonte \`a Baker \cite{baker} et utilisait la 
th\'eorie des surfaces de recouvrement d'Ahlfors.
Des d\'emonstrations plus simples ont \'et\'e 
publi\'ees par Schwick \cite{schwick}, 
Bargmann \cite{bargmann} et Bergweiler \cite{bergweiler}.

On prouvera de fa\c con plut\^ot directe que 
${\cal J}\not=\emptyset $, alors que la d\'emonstration que 
${\cal J}$ est parfait ne sera que une adapation 
d'un argument de \cite{marseille} (voir th.\,1, p.8)
pour les fractions rationnelles,  bas\'e sur le lemme de 
renormalisation de Zalcman (lemme \ref{zalcman}, \cite{zalcman}).
Parmi les nombreux articles \`a ce sujet
le lecteur pourra consulter aussi
\cite{bargmann,bergweiler,bergweiler2,bertelootduval}.

\begin{lemme} 
Soit ${\cal H}:=\{f_{\alpha}\}$ une famille 
de fonctions m\'eromorphes sur un domaine ${\cal W}\subset\CI$:
si ${\cal H}$ n'est pas normale \`a $v\in{\cal W}$, alors il
existe des suites $\{v_n\}\to v$,
$\{r_n\}\downarrow 0$, $\{f_n\}\subset{\cal H}$
et une fonction m\'eromorphe non constante
$h$ sur $\CI$ (\`a d\'eriv\'ee sph\'erique born\'ee sur $\CI$)
telles que $\{f_n(v_n+r_n z)\}\to h$ uniform\'ement sur tout
compact de $ \CI  $.
\QUAN
\labelle{zalcman}
\end{lemme}

Rappelons que 
{\it l'ensemble de Fatou 
$ {\cal F}_f  $}
de
$f$
est l'ensemble des points au 
voisinage desquels
les it\'er\'ees de $f$ 
forment une famille
normale de fonctions holomorphes;
{\it l'ensemble de Julia 
$ {\cal J}_f  $} est le 
compl\'ementaire de 
$ {\cal F}_f  $; $ {\cal F}_f  $ est ouvert alors que 
$ {\cal J}_f  $ est ferm\'e. 

Ces ensembles sont compl\`etement invariants; 
on a aussi $ {\cal J}_f={\cal J}_{f^{\circ n}}  $
et $ {\cal F}_f={\cal F}_{f^{\circ n}}  $ pour tout $n$, voir 
par exemple \cite{bergdyn}, lemme 13.

\begin{lemme}
Soit $f$ une fonction enti\`ere transcendante:
alors il existe $\eta\in\CI$ tel que $f^{\circ 2}(\eta )=\eta $.
\labelle{corps0}
\end{lemme}
{\bf D\'emonstration:} supposons par 
l'absurde qu'il n'existe pas un tel point. 
A fortiori, $f$ n'a
pas de points fixes, ainsi 
$g(z):=({f^{\circ 2}(z)-z})/({f(z)-z})$ 
est une fonction enti\`ere ne prenant
pas la valeur $0$. Gr\^ace au th\'eor\`eme de Picard,
elle prend la valeur $1$: soit 
$z_0\in g^{-1}(1)$. On a $f^{\circ 2}(z_0)=f(z_0)$, ainsi    
$f(z_0)$ est un point fixe de $f$: c'est une contradiction.
\QUAN

\begin{theoreme}
Soit $ f  $ une fonction enti\`ere transcendante:
alors $ {\cal J}_f \not=\emptyset $.
\labelle{corps}
\end{theoreme}
{\sf D\'emonstration} 
Quitte \`a remplacer $f$ par $f^{\circ 2}$, on peut supposer
que $f$ ait au moins un point fixe (lemme \ref{corps0}).\\
{\tt Cas a):} $f$ a au moins
deux points fixes $p_1$ et $p_2$:  
on peut supposer $p_i\in
{\cal F}_f$, $i=1,2$, car sinon il n'y a
rien \`a montrer.
Il existe donc au moins
deux composantes de Fatou
${\cal C}_1$
et ${\cal C}_2$: donc 
$\emptyset\not=
\partial{\cal C}_i\subset
{\cal J}_f $.\\
{\tt Cas b):} $f$ a exactement
un point fixe $p$: supposons par l'absurde 
$ {\cal F}_f =\CI$. 
Or, soit $p$ est un point fixe attractif, soit irrationnellement 
neutre. Dans le premeier cas, $\CI$ est le bassin d'attraction
\`a $p$; dans le deuxi\`eme cas, $\CI$ est un disque de Siegel.
En tout cas, on proive par lin\'earisation de K\"onigs que 
$f$ est injectif sur $\CI$ (voir par exemple \cite{berteloot}, 
th. II.1, II.20, remarque II.21).
C'est une contradiction, car on a suppos\'e $f$ transcendant.
\QUAN
\vskip.2cm
Le th\'eor\`eme suivant est une adaptation du
th\'eor\`eme 1, p.8 de \cite{marseille}
au cas transcendant.

\begin{theoreme}
L'ensemble de Julia $ {\cal J}  $ d'une fonction
enti\`ere transcendante $ f  $ est parfait.
\end{theoreme}
{\sf D\'emonstration} 
Montrons d'abord que $ {\cal J}  $ est infini:
si $ {\cal J}  $ contient au moins deux points,
alors, gr\^ace au thoreme de Picard, 
$ f^{-1}({\cal J} )  $ est infini et on conclut en observant que $  f^{-1}({\cal J} )\subset
{\cal J}     $.
Montrons par ailleurs que l'hypoth\`ese que
 $ {\cal J}  $ contienne exactement un point
$ x  $, elle
m\`ene \`a une contradiction:
en effet, $ f(x)\in{\cal J}  $, donc $ f(x)=x  $.
Analoguement, $ f^{-1}(x)=\{ x  \}  $.
Donc $ z\mapsto (f(z)-x)/(z-x)  $ est une fonction
enti\'ere transcendante ne
prenant pas la valeur $0$, donc
elle prend la valeur $1$ une 
infinit\'e de fois, gr\^ace
au th\'eor\`eme de Picard.
Ainsi
$ f  $ a une infinit\'e de points fixes.
Si un nombre infini de ceux-ci sont contenus
dans $ {\cal J}  $, on a termin\'e; si, au contraire,
il y en a seulement un nombre fini, alors
il existe un nombre infini de composantes de
Fatou, ce qui entra\^\i ne ais\'ement que
$  {\cal J}   $ est \'egalement infini.

Soit maintenant $ \xi\in{\cal J}   $: gr\^ace au lemme \ref{zalcman}, 
il existe des
suites 
$\{v_k\}\to {\xi}$,
$\{r_k\}\downarrow 0$,
$\{f^{\circ n_k}\}
 \subset\{f^{\circ n}\}$
et une fonction  
m\'eromorphe  non constante
$h$
sur $\CI$ 
telles que
$\{f^{\circ n_k}(v_k+r_k z)\}\to h$ 
uniform\'ement
sur tout compact de $ \CI  $;
gr\^ace au lemme de Hurwitz, 
$ h  $ ne prend pas la valeur
$ \infty  $: c'est donc 
une fonction 
enti\`ere.

On peut trouver $ z_0,z_1\in\CI  $ tels que
$ h(z_0)  $ et $ h(z_1)  $ soient deux points
distincts de $ {\cal J}$:
gr\^ace au lemme de Hurwitz, 
il existe deux suites $ \{ z_{0k}  \}\to{z_0}  $
et $ \{ z_{1k}  \}\to{z_1}  $ telles que
$ f^{\circ n_k}(v_k+r_k z_{0k})=h(z_0)  $
et
$ f^{\circ n_k}(v_k+r_k z_{1k})=h(z_1)  $,
ce qui entra\^\i ne 
$v_k+r_k z_{0k}\in{\cal J}$
et
$v_k+r_k z_{1k}\in{\cal J}$ pour tout $k$.
L'une au moins des deux suites 
$\{ v_k+r_k z_{0k}  \}$
et
$\{ v_k+r_k z_{1k}  \}$
n'est pas stationnaire et converge vers
$ \xi  $ dans $ {\cal J}  $.
Cela conclut la d\'emonstration.
\QUAN

\end{document}